\documentclass[11pt]{article}
\usepackage{amsmath,amssymb,amsthm,amsfonts,amstext,amsbsy,amscd}
\usepackage{graphics} \usepackage{graphicx} \usepackage{picins}
\def\<{\langle} \def\>{\rangle}   
   
   \def\Chi{\raise .3ex
  \hbox{\large $\chi$}}   
 \def\({\Bigl (} \def\){\Bigr )}
\newcommand{\be}{\begin{equation}} \newcommand{\ee}{\end{equation}}
\newcommand{\bea}{$$ \begin{array}{lll}} \newcommand{\eea}{\end{array} $$}
\newcommand{\bi}{\begin{itemize}} \newcommand{\ei}{\end{itemize}}

\newtheorem{proposition}{Proposition}

 \renewcommand{\phi}{\varphi}

\begin{document} 
\date{11 janvier 2008} 
\title{{\em A note about conditional Ornstein-Uhlenbeck processes}}
\author{Amel Bentata\footnote{Universit\'e Paris 6, Laboratoire de Probabilit\'es et
    Mod\`eles Al\'eatoires, CNRS-UMR 7599, 16, rue Clisson, 75013 Paris Cedex,
    France}}
\maketitle
\begin{abstract}
  In this short note, the identity in law, which was obtained by P. Salminen \cite{kill}, 
  between on one hand, the Ornstein-Uhlenbeck process with parameter $\gamma$, killed when it reaches $0$, and on the other hand, the  3-dimensional radial 
  Ornstein-Uhlenbeck process killed exponentially at rate $\gamma$ and conditioned to hit $0$, is derived from a simple absolute continuity relationship.
\end{abstract}

\noindent {\it e-mail address:} bentata@clipper.ens.fr \\

\begin{itemize}
  \renewcommand{\labelitemi}{$\bullet$}
\item All probability distributions considered in this note are defined on
  $\mathcal{C}(\mathbb{R}_{+},\mathbb{R})$, where $(X_t,t\geq 0)$ denotes the
  coordinate process, and $\mathcal{F}_t=\sigma\{X_s,s\leq t\}$, its
  natural filtration. 
\item For $\gamma >0$, and $a>0$, we denote by $\mathbb{P}_a^{\gamma}$
  the law of the Ornstein-Uhlenbeck process with parameter $\gamma$,
  starting from $a$, i.e : under $\mathbb{P}_a^{\gamma}$, one has :
  \begin{equation}\label{eq:sde1}
    X_t = a + B_t -\gamma \int_0^t X_s\,ds
  \end{equation}
  for a Brownian motion $(B_t, t\geq 0)$, starting from $0$.

  Note that $(X_t\exp{(\gamma t)}, t\geq 0)$ is a
  $\mathbb{P}_a^{\gamma}$ martingale  which is indeed equal to :
  \begin{equation}\label{eq:mart}
    a+ \int_0^t \exp{(\gamma s)}\,dB_s, \:t \geq 0,
  \end{equation}
  thus leading us to Doob's well-known representation of the
  Ornstein-Uhlenbeck process (see \cite{doob}) :
  \begin{equation}\label{eq:doobrepre}
    X_t=e^{-\gamma t} \Big(a+ \beta \big(\frac{e^{2\gamma t}-1}{2\gamma}\big)\Big),\: t \geq 0,
  \end{equation}
  for another Brownian motion $(\beta(u), u\geq 0)$, starting from $0$.
\item In this note, the 3-dimensional Orstein-Uhlenbeck process $(\vec{X}_t)$ with parameter $\gamma$ and starting from $\vec{a}$, i.e : a solution of (\ref{eq:sde1}) where $(B_t)$ is replaced by a 3-dimensional Brownian motion $(\vec{B}_t)$, and its radial part $R_t=|\vec{X}_t|$ play an important role; $(R_t)$ solves the $SDE$ : 
  \begin{equation}\label{eq:starstar}
    R_t = a + \tilde{B}_t + \int_0^t \frac{ds}{R_s}-\gamma\int_0^t R_s\,ds
  \end{equation} 
  where $a=|\vec{a}|$ and $(\tilde{B}_t, t\geq 0)$ is a 1-dimensional
  Brownian motion, starting from $0$.
\item The main result of this note is the following :
  \begin{proposition}
    Define the probability $\mathbb{Q}_a^{\gamma}$ via :
    \begin{equation}\label{eq:3}
      \mathbb{Q}_a^{\gamma}\vert_{\mathcal{F}_t}=\frac{X_{t\wedge T_o}}{a} \: e^{\gamma t}.\:\mathbb P_a^{\gamma}\vert_{\mathcal{F}_t}
    \end{equation}
    Then, $\mathbb{Q}_a^{\gamma}$ is the law of the 3-dimensional radial
    Ornstein-Uhlenbeck process starting from $a$.
  \end{proposition}
  \begin{proof}
    Due to Girsanov's theorem, and (\ref{eq:mart}), $(B_t)$ considered under $\mathbb{Q}_a^{\gamma}$ solves :
    \begin{equation}\label{eq:croix}
      B_t=\tilde{B}_t + \int_{0}^{t\wedge T_0} \frac{\exp{(\gamma s)}}{X_{s}\,e^{\gamma s}}\:ds
    \end{equation}
    Now, noting that $T_0=\infty$, $\mathbb{Q}_a^{\gamma}$ a.s, we obtain that $X_t$ under $\mathbb{Q}_a^{\gamma}$ satisfies $(\ref{eq:starstar})$, as a consequence of 
    $(\ref{eq:mart})$ and $(\ref{eq:croix})$.
  \end{proof}
\item We now make a few comments :
  \begin{enumerate}
  \item Note that, for $\gamma=0$, (\ref{eq:3}) is nothing else but Doob's
    h-transform relationship between Brownian motion killed when it reaches
    $0$ and the 3-dimensional Bessel process.
  \item We now write $(\ref{eq:3})$ in the equivalent form :
    \begin{equation}\label{eq:seven}
      1_{\{t<T_0\}}\:.\: \mathbb{P}_a^{\gamma}\vert_{\mathcal{F}_t} =
      \frac{a}{X_t}\:e^{-\gamma t}.\:\mathbb{Q}_a^{\gamma}\vert_{\mathcal{F}_t}
    \end{equation}
    which is nothing else but the identity in law obtained by P. Salminen
    (\cite{kill}, Theorem $1$,(i)) between, on the left hand-side, the
    Ornstein-Uhlenbeck process with parameter $\gamma$, killed when it
    reaches $0$, and on the right hand-side the 3-dimensional radial
    Ornstein-Uhlenbeck process killed exponentially at rate $\gamma$, and
    conditioned to hit $0$, that is the density $(\frac{1}{X_t})$
    corresponding to such a conditioning may be seen a posteriori from
    the formula : 
    \begin{equation}
      \mathbb{Q}_a^{\gamma}(F_t
      \frac{1}{X_t})=\mathbb{Q}_a^{\gamma}(\frac{1}{X_t})\,\mathbb{P}_a^{\gamma}(F_t|t<T_0),\: 
      \forall \:F_t \in \mathcal{F}_t\,,
    \end{equation}
    which is a consequence of(\ref{eq:seven})

    It also follows from (\ref{eq:seven}) that $\big(\frac{1}{X_t}\exp{(-\gamma t)}, t\geq 0\big)$ is a strictly local
    martingale with respect to $\mathbb{Q}_a^{\gamma}$, thus extending
    very simply the well-known result for the case $\gamma=0$, when
    $(X_t)$ is a 3-dimensional Bessel process.
    \item Note that (\ref{eq:seven}) leads us to the computation of the semi-group of the killed Ornstein-Uhlenbeck, via :
      \begin{equation}
        \mathbb{E}_a^{\gamma}[\big(\frac{e^{\gamma t}}{a}\big)\:f(X_t)\:1_{\{t\wedge T_0\}}]=\mathbb{Q}_a^{\gamma}[f(X_t)\:\frac{1}{X_t}]\,,
      \end{equation}
      allowing us to recover the result in \cite{fulltext}, p$122$, formula $(53)$.
  \item What about $\gamma<0$ ? Formula (\ref{eq:3}) is still valid thanks to the same arguments; it is closely related
    to Theorem $1$,(ii) of \cite{kill}.
  \item  Related informations about (radial) Ornstein-Uhlenbeck processes, and
    their hitting times, may be found in \cite{Patie}, \cite{strictly}, \cite{Goyor}, \cite{going}.
  \end{enumerate}
\end{itemize}

\bigskip
\noindent{\bf Acknowledgments : } I am very grateful to P. Salminen for a number of discussions during the preparation of this note.

\end{document}